\documentclass [10pt] {article}

 \usepackage [cp1251] {inputenc}
 \usepackage [english] {babel}
 \usepackage {mathrsfs}
 \usepackage {amssymb}
 \usepackage {mathtools}
 \usepackage {amsthm}
 \usepackage {xypic}

 \let \tildeaccent\~   % домашняя страница

 \def \~{\tilde}
 \def \le{\leqslant}
 \def \ge{\geqslant}
 \def \:{\colon}       % функция
 \def \o{\otimes}
 \def \x{\times}
 \def \*{\bullet}
 \def \[{\{}
 \def \]{\}}

 \def \id{{\mathrm{id}}}

 \def \N{\mathbb N}
 \def \Z{\mathbb Z}
 \def \F{\mathbb F}       % простое поле
 \def \k{{\boldsymbol k}} % основное поле

 \def \A{{\mathscr A}} % алгебра Стинрода

 \DeclareMathOperator {\Sq} {Sq}
 \DeclareMathOperator {\St} {St}
 \def \st{\Sigma}
 \DeclareMathOperator {\Hom} {Hom}
 \DeclareMathOperator {\coker} {Coker}
 \DeclareMathOperator {\im} {Im}

 \newcommand* {\head} [1]
 {\subsubsection* {\mathversion{bold}#1}}

 \newcommand* {\subhead} [1]
 {\addvspace\medskipamount
  \noindent {\mathversion{bold}\bf\itshape #1\/}}

 \newenvironment* {claim} [1] []
 {\begin{trivlist}\item [\hskip\labelsep {\bf #1}] \it}
 {\end{trivlist} }

 \newenvironment* {demo} [1] []
 {\begin{trivlist}\item [\hskip\labelsep {\it #1}] }
 {\end{trivlist} }

 \newenvironment* {dem}
 {\begin{trivlist}\item}
 {\end{trivlist} }

\begin {document}

 \title {\large\bf
         Steenrod operations and the diagonal morphism}

 \author {\normalsize\rm
          S.~S.~Podkorytov}

 \date {}

 \maketitle

 \begin {abstract} \noindent
 We show how to find the Steenrod operations in $H^\*(X)$
 (the coefficients in $\mathbb F_p$)
 given
 the diagonal morphism $d_\#:S_\*(X)\to S_\*(X^p)$ and
 the action of the cyclic group $C_p$ on $S_\*(X^p)$.
 Our construction needs no other data
 such as Eilenberg--Zilber morphisms.
 \end {abstract}

 %%%%%%%%%%%%%%%%%%%%%%%%%%%%%%%%%%%%%%%%%%%%%%%%%%%%%%%%%%%%%%

 \head {1. Introduction}

 Fix a prime $p$.
 Chains, cohomology, etc.\ have coefficients in $\F_p$.
 For a (topological) space $X$,
 $S_\*(X)$ is the singular chain complex and
 $d^p_X\:X\to X^p$ is the diagonal map.
 The group $C_p$ of cyclic permutations of the set
 $\{1,\dotsc,p\}$ acts on $X^p$ (on the left)
 in the obvious way.

 \begin {claim} [1.1. Theorem.]
 Let $X$ and $Y$ be spaces.
 Suppose we have a commutative diagram of
 complexes and morphisms
 $$
 \xymatrix {
 S_\*(X)
 \ar[rr]^-{d^p_{X\,\#}}
 \ar[d]_-{f} &&
 S_\*(X^p)
 \ar[d]^-{F} \\
 S_\*(Y)
 \ar[rr]^-{d^p_{Y\,\#}} &&
 S_\*(Y^p),
 }
 $$
 where $F$ preserves the $C_p$-action.
 Then the induced homomorphism $f^*\:H^\*(Y)\to H^\*(X)$
 preserves the action of the Steenrod algebra $\A_p$.
 \end {claim}

 Proof is given in \S~4.
 We do not assume $F$ to be ``$f$ to the power $p$''
 in any sense,
 cf.\ \cite[Satz~4.4]{D}.
 To get an example of a pair $(f,F)$
 satisfying our hypotheses,
 we may take any linear combination of pairs of the form
 $(a_\#,(a^p)_\#)$,
 where
 $a\:X\to Y$ is a continuous map and
 $a^p\:X^p\to Y^p$ is its Cartesian power.
 In this case,
 the assertion of the theorem is obvious.
 Note that
 $f^*$ here need not preserve the cup product.

 Theorem~1.1 implies that
 the Steenrod operations in $H^\*(X)$
 are encoded somehow in the diagonal morphism $d^p_{X\,\#}$
 considered as a morphism from $S_\*(X)$ to some complex of
 $\F_p$-modules with a $C_p$-action.
 We show how to extract Steenrod operations from $d^p_{X\,\#}$,
 see Theorem~4.1.
 Note that
 Steenrod's construction of the operations requires knowing,
 in addition,
 the cross product morphism
 $\xi^p_X\:S_\*(X)^{\o p}\to S_\*(X^p)$
 (see \S~6).

 We conjecture that, conversely,
 the $\A_p$-module $H^\*(X)$ determines
 the diagonal morphism $d^p_{X\,\#}$
 up to coherent quasi-isomorphisms.

 \begin {claim} [Conjecture.]
 Let
 $X$ and $Y$ be spaces
 with finite $\F_p$-homology groups and
 $k\:H^\*(Y)\to H^\*(X)$ be
 an isomorphism of graded $\A_p$-modules.
 Then
 there exists a commutative diagram of
 complexes and morphisms
 $$
 \xymatrix {
 S_\*(X)
 \ar[rr]^-{d^p_{X\,\#}}
 \ar[d]_-{f} &&
 S_\*(X^p)
 \ar[d]^-{F} \\
 S_\*(Y)
 \ar[rr]^-{d^p_{Y\,\#}} &&
 S_\*(Y^p),
 }
 $$
 where
 $F$ preserves the $C_p$-action,
 $f$ and $F$ are quasi-isomorphisms, and
 $f^*=k$.
 \end {claim}

 Most of the paper is taken by known facts about Steenrod's
 construction.
 The only point we state as new,
 besides Theorems 1.1 and 4.1,
 is Corollary~9.2.

 \subhead {Acknoledgements.}
 I thank
 N.~J.~Kuhn and
 S.~V.~Lapin
 for many discussions and
 the referee
 for bringing Hu\`ynh M\`ui's work to my attention and
 other remarks.

 %%%%%%%%%%%%%%%%%%%%%%%%%%%%%%%%%%%%%%%%%%%%%%%%%%%%%%%%%%%%%%

 \head {2. Preliminaries}

 Fix a field $\k$.
 Convention:
 ${\o}={\o_\k}$,
 $\Hom=\Hom_\k$.
 Complexes are over $\k$.
 $\k[n]_\*$ ($n\in\Z$) is the complex with
 $\k[n]_n=\k$ and
 the other terms zero.

 \subhead {Cohomology.}
 For a complex $A_\*$,
 let $h^\*(A_\*)$ be
 the cohomology of the cochain complex $\Hom(A_\*,\k)$.
 We naturally identify $h^n(A_\*)$
 with the $\k$-module of homotopy classes of morphisms
 $A_\*\to\k[n]_\*$.
 For a morphism $f\:A_\*\to\k[n]_\*$,
 we thus have $\[f\]\in h^n(A_\*)$.

 \subhead {Group actions.}
 Let $G$ be a group.
 A $\k G$-module is, clearly, a $\k$-module equipped with a
 ($\k$-linear) $G$-action.
 Given two $\k G$-modules,
 we equip their tensor product (over $\k$) with the diagonal
 $G$-action.
 We equip $\k$ with the trivial $G$-action.

 Let
 $A_\*$ be a $G$-complex
 (i.~e., a complex over $\k G$) and
 $M$ be a $\k G$-module.
 For $n\in\Z$,
 we have the $G$-complex $M[n]_\*$.
 The $G$-homotopy class $\[f\]_G$ of a $G$-morphism
 $f\:A_\*\to M[n]_\*$ is naturally identified with an
 element of the $n$th cohomology $\k$-module of the cochain
 complex $\Hom_G(A_\*,M)$.

 \subhead {Group cohomology.}
 The bar resolution $R_\*(G)$ is a complex of free
 $\k G$-modules.
 $H^\*(G)$ is defined as the cohomology of the cochain complex
 $\Hom_G(R_\*(G),\k)$.

 The augmentation $\epsilon\:R_\*(G)\to\k[0]_\*$ and
 the comultiplication $\Delta\:R_\*(G)\to R_\*(G)\o R_\*(G)$
 are $G$-quasi-isomorphisms.
 $H^\*(G)$ is a graded algebra with
 the unit $1=\[\epsilon\]_G\in H^0(G)$ and 
 the product $H^\*(G)\o H^\*(G)\to H^\*(G)$ induced by
 the composition of morphisms
 $$
 \xymatrix {
 \Hom_G(R_\*(G),\k)\o\Hom_G(R_\*(G),\k)
 \ar[d]^-{(1)} && \\
 \Hom_G(R_\*(G)\o R_\*(G),\k)
 \ar[r]^-{(2)} &
 \Hom_G(R_\*(G),\k),
 }
 $$
 where
 $(1)$ is the obvious product and
 $(2)$ is induced by $\Delta$.

 \subhead {Equivariant cohomology.}
 For
 a $G$-complex $A_\*$ and
 a $\k G$-module $M$,
 let $h_G^\*(A_\*,M)$ be
 the cohomology of the cochain complex
 $\Hom_G(R_\*(G)\o A_\*,M)$.
 If a $G$-morphism $f\:A_\*\to B_\*$ is a quasi-isomorphism,
 then the induced homomorphism
 $h_G^\*(f,\id)\:h_G^\*(B_\*,M)\to h_G^\*(A_\*,M)$ is an
 isomorphism.
 We put $h_G^\*(A_\*)=h_G^\*(A_\*,\k)$.

 $h_G^\*(A_\*)$ is a graded $H^\*(G)$-module with
 the module structure $H^\*(G)\o h_G^\*(A_\*)\to h_G^\*(A_\*)$
 induced by the composition of morphisms
 $$
 \xymatrix {
 \Hom_G(R_\*(G),\k)\o\Hom_G(R_\*(G)\o A_\*,\k)
 \ar[d]^-{(1)} && \\
 \Hom_G(R_\*(G)\o R_\*(G)\o A_\*,\k)
 \ar[r]^-{(2)} &
 \Hom_G(R_\*(G)\o A_\*,\k),
 }
 $$
 where
 $(1)$ is the obvious product and
 $(2)$ is induced by $\Delta$.

 \subhead {The case of a trivial action.}
 Let $A_\*$ be a complex.
 Equip it with the trivial $G$-action.
 The obvious morphism
 $$
 \xymatrix {
 \Hom_G(R_\*(G),\k)\o\Hom(A_\*,\k)
 \ar[r]^-{(1)} &
 \Hom_G(R_\*(G)\o A_\*,\k)
 }
 $$
 induces a graded $\k$-homomorphism
 $$
 \xymatrix {
 H^\*(G)\o h^\*(A_\*)
 \ar[r]^-{(2)} &
 h_G^\*(A_\*),
 }
 $$
 for which we use the notation
 $$
 z\o u\mapsto z\x u.
 $$
 The homomorphism $(2)$ is an $H^\*(G)$-homomorphism.
 If $G$ is finite,
 then $(1)$ and $(2)$ are isomorphisms.

 %%%%%%%%%%%%%%%%%%%%%%%%%%%%%%%%%%%%%%%%%%%%%%%%%%%%%%%%%%%%%%

 \head {3. The algebra $H^\*(C_p)$}

 Here $\k=\F_p$.
 The algebra $H^\*(C_p)$ has a standard basis.
 It consists of classes $e_i\in H^i(C_p)$, $i=0,1,\dotsc$.
 We have $e_ie_j=c_p(i,j)e_{i+j}$,
 where $c_p(i,j)\in\F_p$ is
 $1$ if $pij$ is even and
 $0$ otherwise.
 See \cite[proof of Proposition~9.1, p.~206]{M};
 cf.\ \cite[Ch.~V, \S~5]{SE}.

 %%%%%%%%%%%%%%%%%%%%%%%%%%%%%%%%%%%%%%%%%%%%%%%%%%%%%%%%%%%%%%

 \head {4. Our construction of Steenrod operations}

 Here $\k=\F_p$.
 For each $k\ge0$,
 we define the degree $k$ operation $\st^k\in\A_p$
 by putting:
 $\st^k=\Sq^k$ if $p=2$ and
 $$
 \st^k=
 \begin{cases}
 (-1)^sP^s, &
 k=2(p-1)s, \\
 (-1)^{s+1}\beta P^s, &
 k=2(p-1)s+1, \\
 0 &
 \text{otherwise}
 \end{cases}
 $$
 if $p\ne2$.
 (These $\st^k$ differ from the common $\St^k$ in the sign.)
 For $u\in H^n(X)$,
 we have $\st^ku=0$
 if $n+k>pn$.
 The operations $\st^k$ generate the algebra $\A_p$.

 \begin {claim} [4.1. Theorem.]
 Let
 $X$ be a space and
 $u\in H^n(X)$ ($n\ge0$).
 Consider
 the $C_p$-morphism $d^p_{X\,\#}\:S_\*(X)\to S_\*(X^p)$
 ($C_p$ acts on $S_\*(X)$ trivially)
 and
 the induced homomorphism
 $$
 \xymatrix {
 h_{C_p}^{pn}(S_\*(X)) &&
 h_{C_p}^{pn}(S_\*(X^p))
 \ar[ll]_-{h_{C_p}^{pn}(d^p_{X\,\#})}.
 }
 $$
 Then:
 (a)
 there exists a unique class
 $w\in\im h_{C_p}^{pn}(d^p_{X\,\#})$,
 $$
 w=
 \sum_{i=0}^{pn}
 e_{pn-i}\x v_i,
 \qquad
 v_i\in H^i(X),
 $$
 such that
 $v_i=0$ for $i<n$ and
 $v_n=u$;
 (b)
 we have $v_{n+k}=\st^ku$ for $n\le n+k\le pn$.
 \end {claim}

 Proof is given in \S~9.
 Theorem~4.1 is a way to find/construct the operations $\st^k$.

 \subhead {Proof of Theorem~1.1.}
 Take $u\in H^n(Y)$ ($n\ge0$).
 We show that
 $\st^k f^*(u)=f^*(\st^ku)$ for $n\le n+k\le pn$.

 Consider the commutative diagram
 $$
 \xymatrix {
 h_{C_p}^{pn}(S_\*(X)) &&
 h_{C_p}^{pn}(S_\*(X^p))
 \ar[ll]_-{h_{C_p}^{pn}(d^p_{X\,\#})} \\
 h_{C_p}^{pn}(S_\*(Y))
 \ar[u]^-{h_{C_p}^{pn}(f)} &&
 h_{C_p}^{pn}(S_\*(Y^p))
 \ar[ll]_-{h_{C_p}^{pn}(d^p_{Y\,\#})}
 \ar[u]^-{h_{C_p}^{pn}(F)}.
 }
 $$
 By Theorem 4.1~(a),
 there is a class $w\in\im h_{C_p}^{pn}(d^p_{Y\,\#})$,
 $$
 w=
 \sum_{i=0}^{pn}
 e_{pn-i}\x v_i,
 \qquad
 v_i\in H^i(Y),
 $$
 such that
 $v_i=0$ for $i<n$ and
 $v_n=u$.
 By naturality of the operation $\x$,
 we have
 $$
 h_{C_p}^{pn}(f)(w)=
 \sum_{i=0}^{pn}
 e_{pn-i}\x f^*(v_i).
 $$
 By commutativity of the diagram,
 $h_{C_p}^{pn}(f)(w)\in\im h_{C_p}^{pn}(d^p_{X\,\#})$.
 We have
 $f^*(v_i)=0$ for $i<n$ and
 $f^*(v_n)=f^*(u)$.
 By Theorem 4.1~(b) for $Y$ and $X$,
 $v_{n+k}=\st^ku$ and
 $f^*(v_{n+k})=\st^kf^*(u)$.
 The promised equality follows.
 \qed

 %%%%%%%%%%%%%%%%%%%%%%%%%%%%%%%%%%%%%%%%%%%%%%%%%%%%%%%%%%%%%%

 \head {5. The equivariant power map $\theta_G^r$}

 Let $G$ be a group acting on the set $\{1,\dotsc,r\}$
 ($r\ge0$).
 For $n\in\Z$,
 we define the $G$-module $\k^{(n)}$ as $\k$
 on which an element $g\in G$ acts by multiplication by $s^n$,
 where $s$ is the sign of the permutation corresponding to $g$.

 Let $A_\*$ be a complex.
 $G$ acts on the complex $A_\*^{\o r}$ by
 permuting factors and
 multiplying by $\pm1$ according to the Koszul convention.

 Given a morphism $f\:A_\*\to\k[n]_\*$ ($n\in\Z$),
 we introduce the $G$-morphism
 $$
 \xymatrix {
 f_G^r\:
 R_\*(G)\o A_\*^{\o r}
 \ar[rr]^-{\epsilon\o f^{\o r}} &&
 \k[0]_\*\o\k[n]_\*^{\o r}
 \ar[r]^-{(1)} &
 \k^{(n)}[rn]_\*,
 }
 $$
 where $(1)$ is the $G$-isomorphism
 given by $1\o1^{\o r}\mapsto1$.
 We have
 $\[f\]\in h^n(A_\*)$ and
 $\[f_G^r\]_G\in h_G^{rn}(A_\*^{\o r},\k^{(n)})$.

 \begin {claim} [5.1. Lemma.]
 Let $f_0,f_1\:A_\*\to\k[n]_\*$ be homotopic morphisms.
 Then the $G$-morphisms $(f_i)_G^r$, $i=0,1$, are
 $G$-homotopic.
 \end {claim}

 This follows from \cite[Lemma 1.1~(ii)]{M};
 cf.\ \cite[Ch.~VII, Lemma~2.2]{SE}.

 \begin {demo} [Proof.]
 Let $L_\*$ be the complex with
 a basis consisting of $b\in L_{n-1}$ and $c_0,c_1\in L_n$ and
 the differential given by $\partial c_i=(-1)^ib$.
 Define morphisms
 $s\:\k[n]_\*\to L_\*$ and
 $t_i\:L_\*\to\k[n]_\*$, $i=0,1$,
 by
 $s(1)=c_0+c_1$ and
 $t_i(c_j)=\delta_{ij}$ (the Kronecker delta).
 We have $t_i\circ s=\id$.
 Thus,
 for each $i$,
 we have the commutative diagram
 $$
 \xymatrix {
 &
 h_G^{rn}(L_\*^{\o r},\k^{(n)})
 \ar[dl]_-{s^\#}
 & \\
 h_G^{rn}(\k[n]_\*^{\o r},\k^{(n)}) & &
 h_G^{rn}(\k[n]_\*^{\o r},\k^{(n)})
 \ar[ll]_-{\id}
 \ar[ul]_-{t_i^\#}
 }
 $$
 (we put $s^\#=h_G^{rn}(s^{\o r},\id)$, etc.).
 Since $s$ is a quasi-isomorphism,
 $s^\#$ is an isomorphism.
 It follows that $t_0^\#=t_1^\#$.
 Put $e=\id\:\k[n]_\*\to\k[n]_\*$.
 We have $\[(t_i)_G^r\]_G=t_i^\#(\[e_G^r\]_G)$.
 Thus $\[(t_0)_G^r\]_G=\[(t_1)_G^r\]_G$.

 It follows from the definition of chain homotopy that
 there exists a morphism $F\:A_\*\to L_\*$ such that
 $f_i=t_i\circ F$, $i=0,1$.
 Consider the homomorphism $F^\#\:
 h_G^{rn}(L_\*^{\o r},\k^{(n)})\to
 h_G^{rn}(A_\*^{\o r},\k^{(n)})$.
 We have $\[(f_i)_G^r\]_G=F^\#(\[(t_i)_G^r\]_G)$.
 Thus $\[(f_0)_G^r\]_G=\[(f_1)_G^r\]_G$.
 \qed
 \end {demo}

 Lemma~5.1 allows us to introduce the map
 $$
 \theta_G^r\:
 h^n(A_\*)\to
 h_G^{rn}(A_\*^{\o r},\k^{(n)}),
 \quad
 \[f\]\mapsto(-1)^{nr(r-1)/2}\[f_G^r\]_G.
 $$
 It is not additive in general.

 %%%%%%%%%%%%%%%%%%%%%%%%%%%%%%%%%%%%%%%%%%%%%%%%%%%%%%%%%%%%%%

 \head {6. Steenrod's construction of $\st^k$}

 Here $\k=\F_p$.
 Define the {\it Steenrod numbers\/} $a_p(n)\in\F_p$, $n\in\N$:
 put
 $a_p(n)=1$ if $p=2$ and
 $a_p(n)=(-1)^{qn(n+1)/2}(q!)^n$ if $p=2q+1$.
 We have $a_p(n)\ne0$.

 Let $X$ be a space.
 We have the cross product quasi-isomorphism
 $\xi^p_X\:S_\*(X)^{\o p}\to S_\*(X^p)$,
 which preserves the action of $C_p$
 (in fact, of the whole symmetric group $\Sigma_p$).
 Consider the diagram of $C_p$-complexes
 $$
 \xymatrix {
 S_\*(X)
 \ar[rr]^-{d^p_{X\,\#}} &&
 S_\*(X^p) &&
 S_\*(X)^{\o p},
 \ar[ll]_-{\xi^p_X}
 }
 $$
 where the $C_p$-action in $S_\*(X)$ is trivial.
 Given a class $u\in H^n(X)=h^n(S_\*(X))$ ($n\ge0$),
 we have the class
 $\theta_{C_p}^p(u)\in h_{C_p}^{pn}(S_\*(X)^{\o p})$
 (note that $\F_p^{(n)}=\F_p$ as an $\F_pC_p$-module).
 Define the classes $\Psi(u)$ and $\Phi(u)$:
 $$
 \xymatrix {
 \Phi(u) &&
 \Psi(u)
 \ar@{|->}[ll]
 \ar@{|->}[rr] &&
 a_p(n)^{-1}\theta_{C_p}^p(u) \\
 h_{C_p}^{pn}(S_\*(X)) &&
 h_{C_p}^{pn}(S_\*(X^p))
 \ar[ll]_-{h_{C_p}^{pn}(d^p_{X\,\#})}
 \ar[rr]^-{h_{C_p}^{pn}(\xi^p_X)} &&
 h_{C_p}^{pn}(S_\*(X)^{\o p}).
 }
 $$
 They are well-defined
 since $h_{C_p}^{pn}(\xi^p_X)$ is an isomorphism.

 We have
 $$
 \Phi(u)=
 \sum_{i=0}^{pn}
 e_{pn-i}\x\phi_i(u)
 $$
 for some $\phi_i(u)\in h^i(S_\*(X))=H^i(X)$, $i=0,\dotsc,pn$.

 \begin {claim} [6.1. Fact.]
 We have
 $\phi_i(u)=0$ for $i<n$,
 $\phi_n(u)=u$, and
 $\phi_{n+k}(u)=\st^ku$ for $n\le n+k\le pn$.
 \end {claim}

 \begin {dem}
 This is the construction of Steenrod operations given in
 \cite[\S\S\ 2, 5, 7, 8]{M},
 note
 Remarks~7.2,
 Theorem 7.9~(i),
 Propositions 8.1 and 2.3~(iv),
 and also formula (2) in the proof of Proposition~9.1 there.
 We used unnormalized chains,
 but we could use normalized ones equivalently,
 thus following \cite{M}.
 Cf.\ also \cite[Ch.~VII]{SE}.
 \qed
 \end {dem}

 %%%%%%%%%%%%%%%%%%%%%%%%%%%%%%%%%%%%%%%%%%%%%%%%%%%%%%%%%%%%%%

 \head {7. The transfer and the functor $\~h_G^\*$}

 Let $G$ be a finite group.
 Put
 $$
 M=\sum_{g\in G}g\in\k G.
 $$
 Let $A_\*$ be a $G$-complex.
 The morphism
 $\epsilon\o M\:R_\*(G)\o A_\*\to\k[0]_\*\o A_\*=A_\*$ induces
 a morphism $\Hom(A_\*,\k)\to\Hom_G(R_\*(G)\o A_\*,\k)$.
 The induced $\k$-homomorphism on the cohomology
 $t\:h^\*(A_\*)\to h_G^\*(A_\*)$ is called the {\it transfer\/}
 (\cite[Ch.~V, \S~7]{SE}, cf.\ \cite[Ch.~III, \S~9]{B}).

 Equip $h^\*(A_\*)$ with a (trivial) $H^\*(G)$-module structure
 using the $\k$-algebra homomorphism $H^\*(G)\to H^\*(1)=\k$
 induced by the group inclusion $1\to G$.
 Then $t$ becomes an $H^\*(G)$-homomorphism
 (cf.\ \cite[Ch.~V, (3.8)]{B}).

 Put $\~h_G^\*(A_\*)=\coker t$.
 This is an $H^\*(G)$-module.
 Obviously,
 $\~h_G^\*$ is a functor.
 It respects quasi-isomorphisms.

 The following lemmas are immediate.

 \begin {claim} [7.1. Lemma.]
 Let $A_\*$ be a complex.
 Consider the $G$-complex $B_\*=A_\*\o\k G$.
 Then $\~h_G^\*(B_\*)=0$.
 \qed
 \end {claim}

 \begin {claim} [7.2. Lemma.]
 Let $A_\*$ be a complex.
 Equip it with the trivial $G$-action.
 Then the transfer $t\:h^\*(A_\*)\to h_G^\*(A_\*)$ is given by
 $t(u)=\lvert G\rvert\x u$.
 \qed
 \end {claim}

 %%%%%%%%%%%%%%%%%%%%%%%%%%%%%%%%%%%%%%%%%%%%%%%%%%%%%%%%%%%%%%

 \head {8. The functor $\~h_{C_p}^\*(?^{\o p})$ and
        the map $\~\theta_{C_p}^p$}

 Here we follow \cite[Ch.~II, proof of Theorem~3.7]{H}.

 We have a prime $p$.
 We show that
 the functor $\~h_{C_p}^m(?^{\o p})$ takes sums to products.
 It follows that
 it is additive.

 \begin {claim} [8.1. Lemma.]
 Let $A_{i\,\*}$, $i\in I$, be a family of complexes.
 Put
 $$
 A_\*=\bigoplus_{i\in I}A_{i\,\*}.
 $$
 Let $k_i\:A_{i\,\*}\to A_\*$ be the canonical morphisms.
 Then
 $$
 \xymatrix {
 \~h_{C_p}^m(A_\*^{\o p})
 \ar[rrr]^-{(\~h_{C_p}^m(k_i^{\o p}))_{i\in I}} &&&
 \prod_{i\in I}\~h_{C_p}^m(A_{i\,\*}^{\o p})
 }
 $$
 is an isomorphism for each $m\in\Z$.
 \end {claim}

 \begin {demo} [Proof.]
 Put $J=I^p\setminus\im d$,
 where $d\:I\to I^p$ is the diagonal map.
 $C_p$ acts on $J$ by permuting coordinates.
 Consider the complex
 $$
 N_\*=
 \bigoplus_{(i_1,\dotsc,i_p)\in J}
 A_{i_1\,\*}\o\dotso\o A_{i_p\,\*}.
 $$
 Equip it with the obvious $C_p$-action.
 We have the $C_p$-isomorphism
 $$
 \xymatrix {
 A_\*^{\o p} &&&
 \big(
 \bigoplus_{i\in I}
 A_{i\,\*}^{\o p}
 \big)
 \oplus N_\*,
 \ar[lll]_-{((k_i^{\o p})_{i\in I},l)}
 }
 $$
 where $l=(k_{i_1}\o\dotso\o k_{i_p})_{(i_1,\dotsc,i_p)\in J}\:
 N_\*\to A_\*^{\o p}$.
 Obviously,
 the functor $\~h_{C_p}^m$ takes sums to products.
 Thus
 $$
 \xymatrix {
 \~h_{C_p}^m(A_\*^{\o p})
 \ar[rrrr]^-{((\~h_{C_p}^m(k_i^{\o p}))_{i\in I},
              \~h_{C_p}^m(l))}
 &&&&
 \big(
 \prod_{i\in I}
 \~h_{C_p}^m(A_{i\,\*}^{\o p})
 \big)
 \oplus\~h_{C_p}^m(N_\*)
 }
 $$
 is an isomorphism.
 It remains to show that $\~h_{C_p}^m(N_\*)=0$.

 Choose a section $s\:J/C_p\to J$ and
 put $J'=\im s\subseteq J$.
 Consider the complex
 $$
 N'_\*=
 \bigoplus_{(i_1,\dotsc,i_p)\in J'}
 A_{i_1\,\*}\o\dotso\o A_{i_p\,\*}.
 $$
 Since $p$ is prime,
 $C_p$ acts on $J$ freely.
 It follows that $N_\*\cong N'_\*\o\k C_p$.
 By Lemma~7.1,
 $\~h_{C_p}^m(N_\*)=0$.
 \qed
 \end {demo}

 Set $\k=\F_p$.
 Given
 a complex $A_\*$ and
 a class $u\in h^n(A_\*)$ ($n\in\Z$),
 we define $\~\theta_{C_p}^p(u)\in\~h_{C_p}^p(A_\*^{\o p})$
 as the image of $\theta_{C_p}^p(u)\in h_{C_p}^p(A_\*^{\o p})$
 under the projection
 (note that $\F_p^{(n)}=\F_p$ as an $\F_pC_p$-module).

 \begin {claim} [8.2. Lemma.]
 Let $A_\*$ be a complex.
 Then:
 (a)
 the map
 $\~\theta_{C_p}^p\:h^n(A_\*)\to\~h_{C_p}^{pn}(A_\*^{\o p})$ is
 $\F_p$-linear;
 (b)
 the $H^\*(C_p)$-homomorphism
 $$
 H^\*(C_p)\otimes h^\*(A_\*)\to
 \~h_{C_p}^\*(A_\*^{\o p}),
 \quad
 1\otimes u\mapsto
 \~\theta_{C_p}^p(u),
 \
 u\in h^n(A_\*),
 \
 n\in\Z,
 $$
 is a (non-graded) isomorphism.
 \end {claim}

 \begin {demo} [Proof.]
 (a)
 The map
 $\~\theta_{C_p}^p\:h^n(A_\*)\to\~h_{C_p}^{pn}(A_\*^{\o p})$ is
 natural in $A_\*$.
 Its source and target are additive in $A_\*$
 (obvious for $h^n(A_\*)$,
 Lemma~8.1 for the $\~h_{C_p}^{pn}(A_\*^{\o p})$).
 It follows that it is additive and thus $\F_p$-linear.

 (b)
 There exists a quasi-isomorphism
 $$
 \xymatrix {
 A_\* &&
 \bigoplus_{i\in I}E_{i\,\*},
 \ar[ll]_-{q=(q_i)_{i\in I}}
 }
 $$
 where $E_{i\,\*}=\F_p[m_i]_\*$, $m_i\in\Z$.
 We have the commutative diagram
 $$
 \xymatrix {
 h^n(A_\*)
 \ar[rr]^-{(h^n(q_i))_{i\in I}}
 \ar[d]_-{\~\theta_{C_p}^p} &&
 \prod_{i\in I}h^n(E_{i\,\*})
 \ar[d]^-{\prod_{i\in I}
          (\~\theta_{C_p}^p\:
           h^n(E_{i\*})\to
           \~h_{C_p}^{pn}(E_{i\*}^{\o p}))} \\
 \~h_{C_p}^{pn}(A_\*^{\o p})
 \ar[rr]^-{(\~h_{C_p}^{pn}(q_i^{\o p}))_{i\in I}} &&
 \prod_{i\in I}
 \~h_{C_p}^{pn}(E_{i\,\*}^{\o p}).
 }
 $$
 The functors $h^n$ and $\~h_{C_p}^{pn}(?^{\o p})$
 respect quasi-isomorphisms
 (obvious)
 and take sums to products
 (obvious for $h^n$,
 Lemma~8.1 for $\~h_{C_p}^{pn}(?^{\o p})$).
 Thus the horizontal arrows are isomorphisms.
 Tensoring by $H^\*(C_p)$ commutes with products
 since $H^k(C_p)$ are finite.
 Thus all reduces to the case $A_\*=\F_p[m]_\*$,
 where the assertion is verified immediately.
 \qed
 \end {demo}

 %%%%%%%%%%%%%%%%%%%%%%%%%%%%%%%%%%%%%%%%%%%%%%%%%%%%%%%%%%%%%%

 \head {9. Back to topology}

 Here $\k=\F_p$.
 Let $X$ be a space.

 \begin {claim} [9.1. Lemma.]
 (a)
 The maps $\Phi\:H^n(X)\to h_{C_p}^{pn}(S_\*(X))$, $n\ge0$, are
 $\F_p$-linear.
 (b)
 The image of the $H^\*(C_p)$-homomorphism
 $$
 \xymatrix {
 h_{C_p}^\*(S_\*(X)) &&
 h_{C_p}^\*(S_\*(X^p))
 \ar[ll]_-{h_{C_p}^\*(d_{X\,\#}^p)}
 }
 $$
 is $H^\*(C_p)$-generated by classes of the form $\Phi(u)$,
 $u\in H^n(X)$, $n\ge0$.
 \end {claim}

 This is known,
 see \cite[Ch.~II, proof of Theorem~3.7]{H}.
 The assertion (a) follows also from
 Fact~6.1 and
 the well-known linearity of $\st^k$.
 Our proof follows \cite{H}.

 \begin {demo} [Proof.]
 We have the commutative diagram
 $$
 \xymatrix {
 \Phi(u) &&
 \Psi(u)
 \ar@{|->}[ll]
 \ar@{|->}[rr] &&
 a_p(n)^{-1}\theta_{C_p}^p(u) \\
 h_{C_p}^{pn}(S_\*(X))
 \ar[d]^-{r} &&
 h_{C_p}^{pn}(S_\*(X^p))
 \ar[ll]_-{h_{C_p}^{pn}(d^p_{X\,\#})}
 \ar[rr]^-{h_{C_p}^{pn}(\xi^p_X)}
 \ar[d] &&
 h_{C_p}^{pn}(S_\*(X)^{\o p})
 \ar[d] \\
 \~h_{C_p}^{pn}(S_\*(X)) &&
 \~h_{C_p}^{pn}(S_\*(X^p))
 \ar[ll]_-{\~h_{C_p}^{pn}(d^p_{X\,\#})}
 \ar[rr]^-{\~h_{C_p}^{pn}(\xi^p_X)} &&
 \~h_{C_p}^{pn}(S_\*(X)^{\o p}),
 }
 $$
 where the vertical arrows are the projections.
 By Lemma~7.2,
 the transfer for $S_\*(X)$ is zero
 (cf.\ \cite[Ch.~VII, proof of Lemma~4.1]{SE}).
 Thus $r$ is an isomorphism.
 Since $\xi^p_X$ is a $C_p$-quasi-isomorphism,
 $\~h_{C_p}^\*(\xi^p_X)$ is an isomorphism.
 By Lemma~8.2,
 the maps
 $\~\theta_{C_p}^p\:H^n(X)\to\~h_{C_p}^{pn}(S_\*(X)^{\o p})$ are
 $\F_p$-linear and
 the $H^\*(C_p)$-module $\~h_{C_p}^\*(S_\*(X)^{\o p})$ is
 generated by classes of the form $\~\theta_{C_p}^p(u)$,
 $u\in H^n(X)$, $n\ge0$.
 The desired assertions follow.
 \qed
 \end {demo}

 \begin {claim} [9.2. Corollary.]
 Take a class $w\in\im h_{C_p}^m(d_{X\,\#}^p)$ ($m\ge0$),
 $$
 w=
 \sum_{i=0}^m
 e_{m-i}\x v_i,
 \qquad
 v_i\in H^i(X).
 $$
 If $v_i=0$ for $i\le m/p$,
 then $w=0$.
 \end {claim}

 \begin {demo} [Proof.]
 By Lemma~9.1 (a, b),
 we have
 \begin{equation}
 w=
 \sum_{k=0}^{[m/p]}
 e_{m-pk}\Phi(u_k)
 \tag{$*$}
 \end{equation}
 for some $u_k\in H^k(X)$.
 Using the formulas
 $$
 \Phi(u_k)=
 \sum_{i=0}^{pk}
 e_{pk-i}\x\phi_i(u_k)
 $$
 and
 $e_ie_j=c_p(i,j)e_{i+j}$,
 we get
 $$
 v_i=
 \sum_{k:\,i\le pk\le m}
 c_p(m-pk,pk-i)\phi_i(u_k).
 $$
 We
 take successively $n=0,\dotsc,[m/p]$ and
 show that $u_n=0$.
 On each step,
 $u_k=0$ for $k<n$.
 By Fact~6.1,
 $\phi_n(u_n)=u_n$ and
 $\phi_n(u_k)=0$ for $k>n$.
 Therefore, $v_n=c_p(m-pn,pn-n)u_n$.
 We have $c_p(m-pn,pn-n)=1$
 since $p(pn-n)$ is even.
 Thus $v_n=u_n$.
 By assumption,
 $v_n=0$.
 Thus $u_n=0$.
 \qed
 \end {demo}

 A similar reasoning shows that
 the presentation ($*$) is unique for any
 $w\in\im h_{C_p}^m(d_{X\,\#}^p)$.
 Thus we get a (non-graded) isomorphism
 $\im h_{C_p}^\*(d_{X\,\#}^p)\cong H^\*(C_p)\otimes H^\*(X)$.
 It follows also that $\~h_{C_p}^\*(d_{X\,\#}^p)$ is injective.
 This stuff is known,
 see \cite[Ch.~II, proof of Theorem~3.7]{H}.

 \subhead {Proof of Theorem~4.1.}
 The uniqueness follows from Corollary~9.2.
 By Fact~6.1,
 $w=\Phi(u)$ satisfies the conditions of the parts (a) and (b).
 Thus we are done.
 \qed

 %%%%%%%%%%%%%%%%%%%%%%%%%%%%%%%%%%%%%%%%%%%%%%%%%%%%%%%%%%%%%%

 %% %% %% %% %% %% %% %% %% %% %% %% %% %% %% %% %% %% %% %% %%

 \medskip

 {\noindent \tt ssp@pdmi.ras.ru}

 {\noindent \tt http://www.pdmi.ras.ru/\tildeaccent{}ssp}

 \end {document}